# The distribution of prime numbers and tuples in the natural numbers

VICTOR VOLFSON

ABSTRACT. We studied two probabilistic models of the distribution of primes in the natural number [1]. The paper considers the third probabilistic model of the distribution of primes in the natural number. The author proved that the results obtained by considering three probabilistic models are similar. It is shown that the accuracy of these probabilistic models exceeds the Riemann hypothesis with a certain probability. The author constructed probabilistic models for the analysis of the average distance between neighboring prime numbers and k-tuples and found probability distributions of these random variables. The estimates of these random variables are obtained based on these probabilistic models. The resulting estimates tested on a large volume of statistical data.

1. INTRODUCTION

We studied two probabilistic models of the distribution of primes in the natural number [1]. We consider the third probabilistic model of the distribution of prime numbers in this paper. It will be shown that the results obtained using these probabilistic models are very similar at large intervals of natural numbers.

Probabilistic methods may be advantageously used also in the analysis of the distance between adjacent primes. It is proved [2] that there are over $C_1 \ln(x)$ various increments $d_i = p_{i+1} - p_i < C_2 \ln(x), (d_i \ll \ln(x))$ (for suitable constants $C_1, C_2$), formed by prime numbers $p_{i+1}, p_i$ (not exceeding value $x$). Consequently, it was shown that the increase $d_i \ll \ln(x)$ give a finite fraction of all the increments between primes i.e. the probability of this event $A$ is positive $P(A) > 0$. Thus, the probability that $d_i \gg \ln(x)$ is $1 - P(A)$.

Probabilistic techniques have been used to find the maximum distance between adjacent primes not exceeding value $x$ [3], and it was shown, that the inequality $sup_i(d_i) \leq C_3 \ln^2(x)$. The probabilistic methods are used to analyze the average distance between neighboring primes in the paper.

---

Probabilistic methods used to find the maximum distance between prime $k$ - tuples [4], and constructed the model to estimate the probability distribution of prime $k$ -tuples in natural numbers [5]. The probabilistic methods are used to find the average distance between adjacent $k$ - tuples in the paper.

The ideas of the probabilistic methods are well explained Kramer in his work [6]: «It is often possible the following interesting application of probabilistic reasoning in the arguments, related to the asymptotic properties of arithmetic functions. If, for example, we are interested in the distribution of the sequence of $S$ integers, then we regarded $S$ as an element of an infinite class $C$ of sequences, which can be specifically interpreted as the possible implementation of some of the game case. Then, we can prove (in many cases), that with probability equal to 1, a ratio $R$ is performed in $C$ or in the precise mathematical sense of "almost all" of the sequence from $C$ are satisfied with relationship - $R$. Of course, we cannot, in general, conclude that the relationship $R$ is performed for each sequence $S$, but the results predicted by this method can sometimes be rigorously proved by other methods».

Thus, following Cramer's approach, we can assume that the difference in values between the actual number of primes not exceeding $x$, and calculated by the formula $Li(x)$ has a normal distribution (based on the first and second probabilistic model of the distribution of prime numbers in the natural number [1]).

Also, based on the probability distribution model prime $k$ - tuples in the natural numbers [5], we can assume that the difference between the actual number of prime $k$ - tuples in the natural numbers (not exceeding $x$) and calculated (on the basis of the of Hardy-Littlewood conjecture) has a normal distribution.

2. PROBABALISTIC MODELS OF THE DISTRIBUTION OF PRIMES IN THE NATURAL NUMBERS

It is assumed [1] (in the first probabilistic model) that somebody at random chooses a ball from the basket containing balls with natural numbers from 1 to $x$. If the number of the ball is a prime, then the random variable is assigned the value - 1, if not- 0. The selected ball again returns to the basket in this probabilistic model. Therefore, there is a possibility to choose one and the same ball a few times in the model. This situation does not happen in a real time, where calculated the number of primes in the range of the natural numbers from 1 to $x$.

Let us consider the probabilistic model, which has no such drawbacks - with choice without returning the ball. We do not get the binomial distribution of the sum of random variables as in the first probabilistic model. We get the hypergeometric distribution in this choice.

We consider the urn containing $M$ different balls. There are $M_1$ white balls among them and $M_2$ - black ($M_1 + M_2 = M$). Suppose that somebody chooses $n$ balls from the urn ($n < M$) without returning balls back, then the probability of events that somebody choose $n_1$ white balls and $n_2$ black balls is equal to:

$$P(B_{n_1,n_2}) = C_{M_1}^{n_1} \cdot C_{M_2}^{n_2} / C_M^n .$$

It can be shown [7] that if $M$ and $M_1$ tend to infinity and if the limit $M_1/M$ tends to probability- $p$ (therefore $M_2/M$ tends to $1-p$), then the limit $P(B_{n_1,n_2})$ tends to $C_{n_1+n_2}^{n_1} p^{n_1}(1-p)^{n_2}$.

Thus, (under these assumptions) the hypergeometric distribution tends to the binomial distribution. The final choice without returning must give almost the same result as the choice with returning for large values $M, M_1$.

Let us return back to the probabilistic model. We take the first $M$ natural numbers. Let assume that $M_1$ numbers are primes among them. The value $M_1$ tends to infinity when the value $M$ tends to infinity too. The ratio $M_1/M$ is equal to $p = Li(M)/M$ (the probability of randomly choose positive integer, which does not exceed value $M$, to be prime) for large value $M$.

Thus, the first probabilistic model gives almost the same effect as discussed probabilistic model without returning for large values $M$.

The second probabilistic model [1] considers the random variables - similar number of primes not exceeding a value $x$.

It was shown (in the second probabilistic model) that the random variable has a normal distribution with mean equal to $a = Li(x)$ and standard deviation equal to:

$$\sigma_1 = \sqrt{Li(x) - Li_2(x)} , \qquad (2.1)$$

where $Li_k(x) = \int_2^x dt / \ln^k(t)$.

It was shown (in the first probabilistic model) for the random variable with the same expectation, that the standard deviation is equal to:

$$\sigma = \sqrt{Li(x) - Li^2(x)/x}. \tag{2.2}$$

Based on (2.1) and (2.2) we prove the following statements:

1. $\sigma_1^2 \sim \sigma^2$. $\qquad(2.3)$

2. $\sigma^2 - \sigma_1^2 = x/\ln^3(x) + O(x/\ln^4(x))$. $\qquad(2.4)$

We prove (2.3) the first.

After integration by parts:

$$Li(x) = x/\ln(x) + Li_2(x). \tag{2.5}$$

Based on (2.5):

$$Li^2(x)/x = x/\ln^2(x) + 2 \cdot Li_2(x)/\ln(x) + Li_2^2(x)/x \tag{2.6}$$

After integration by parts:

$$1/\ln(x)\; Li_2(x) = x/\ln^2(x) + Li_3(x). \tag{2.7}$$

Substituting (2.7) into (2.6):

$$Li^2(x)/x = x/\ln^2(x) + 2x/\ln^3(x) + O(x/\ln^4(x)). \tag{2.8}$$

Based on $Li_3(x) = x/\ln^3(x) + O(x/\ln^4(x))$ and (2.7) we obtain:

$$Li_2(x) = x/\ln^2(x) + x/\ln^3(x) + O(x/\ln^4(x)). \tag{2.9}$$

Based on (2.5), (2.8) and (2.9) we obtain:

$$\lim_{x \to \infty} \sigma^2/\sigma_1^2 = \lim_{x \to \infty} \frac{Li(x) - Li^2(x)/x}{Li(x) - Li_2(x)} = \lim_{x \to \infty} \frac{x/\ln(x) - x/\ln^3(x) + O(x/\ln^4(x))}{x/\ln(x) + O(x/\ln^4(x))} = 1$$

We prove (2.4) now. Based on (2.8) and (2.9) we obtain:

$$\sigma^2 - \sigma_1^2 = Li^2(x)/x - Li_2(x) =$$

$$x/\ln^2(x)+2x/\ln^3(x)-x/\ln^2(x)-x/\ln^3(x)+O(x/\ln^4(x))=x/\ln^3(x)+O(x/\ln^4(x))$$

Assertions (2.2) and (2.3) are in agreement with the calculated data in Table 1.

Let compare the following measures of accuracy of these models and the Riemann conjecture (Table 1): the number of primes that do not exceed a value $x$ - $\pi(x)$, the whole value of the difference of values -$[Li(x)-\pi(x)]$, the whole value of the standard deviation for the first probabilistic model -$[\sqrt{Li(x)-Li^2(x)/x}]$, the whole value of the standard deviation for the second probabilistic model -$[\sqrt{Li(x)-\int_2^x dt/\ln^2(t)}]$, the whole the value of the maximum deviation of Riemann conjecture -$[\sqrt{x}\ln(x)/8\pi]$.

Таблица 1

| x | The number of prime numbers not exceeding value x | The whole value of the difference $[Li(x)-\pi(x)]$ | The whole value of the standard deviation for the first probabilistic model | The whole value of the standard deviation for the second probabilistic model | The whole value of the maximum deviation $\pi(x)$ from $Li(x)$ of Riemann conjecture |
|---|---|---|---|---|---|
| $10^8$ | 5761455 | 754 | 2330 | 2329 | 7333 |
| $10^9$ | 50847534 | 1701 | 7091 | 7089 | 26087 |
| $10^{10}$ | 455052511 | 3104 | 20841 | 20839 | 91663 |
| $10^{11}$ | 4118054813 | 11588 | 62836 | 62834 | 318851 |
| $10^{12}$ | 37607912018 | 38263 | 190246 | 190239 | 1099961 |

The conclusions from Table 1:

1. The whole value of the difference $[Li(x) - \pi(x)]$ is relatively small and fit into a single value of the standard deviation of the first and second probabilistic models.

2. Values of standard deviation for the first and second (model Kramer) probabilistic models (for large values $x$) are almost equal, which corresponds to the statements (2.2), (2.3). However, Kramer make a heuristic assume in his model, that the probability of a large value $x$ to be prime is equal to $1/\ln(x)$, but in the first probabilistic model used the proving assertion that the probability to choose a prime from the natural numbers (not exceeding $x$) is equal to $Li(x)/x$.

3. The value of the maximum deviation $\pi(x)$ from $Li(x)$ in the Riemann conjecture is significantly higher than the values for the standard deviations of the first and second probabilistic models.

These conclusions are confirmed in [8] on the basis of the analysis of zeros of the Riemann zeta function (for small v the distribution $[Li(x) - \pi(x)]$ is close to a normal distribution).

The probability density distribution of the random variable $x_1$ - analogue of the number of primes that do not exceed the value $x$, $(x \geq x_0, x_0 -$ *the number of Schiusa*) determinates by the formula of the normal distribution (for the first probabilistic model):

$$P(x_1) = \frac{1}{\sigma\sqrt{2\pi}} \cdot e^{\frac{-(x_1-a)^2}{2(\sigma)^2}}, \qquad (2.10)$$

where $a = Li(x), \sigma = \sqrt{Li(x) - Li(x)^2/x}$.

Similarly, (2.10) is the probability density distribution of the random variable $x_2$ - analogue of the number of primes that do not exceed the value $x$ $(x < x_0)$, determines by the formula:

$$P(x_2) = \frac{2}{\sigma\sqrt{2\pi}} \cdot e^{\frac{-(x_2-a)^2}{2(\sigma)^2}} \text{ for } x_2 \leq a \text{ and } P(x_2) = 0 \text{ for } x_2 > a. \qquad (2.11)$$

Note that the inequality holds for large values $x$:

$$x/\ln(x) < Li(x) - 3\sigma < ... < Li(x). \qquad (2.12)$$

Based on (2.11) the function of the probability distribution of the random variable $x_2$ has the form:

$$F(x_2) = \frac{2}{\sigma\sqrt{2\pi}} \cdot \int_{-\infty}^{x_2} e^{\frac{-(t-a)^2}{2(\sigma)^2}} dt \text{ for } x_2 \leq a \text{ and } F(x_2) = 1 \text{ for } x_2 > a. \qquad (2.13)$$

Let us find the distribution density of the random variable $Y = X_2/x$, where $x < x_0$ ($x_0$ - the number of Schiusa). We use the formula:

$$P_Y(Y) = P_{X_2}[g(Y)]g'(Y) \qquad (2.14)$$

(where $g(Y)$ - the inverse function) to determine the probability density function of the random variable [7]:

In this case:

$$g(Z) = 1/Z, g'(Z) = -1/Z^2, g(Y) = x \cdot Y, g'(Y) = x.$$

Therefore, based on (2.13), (2.14) if $X_2 \leq a$ or $Y \leq a/x = Li(x)/x$ we obtain:

$$P_Y(Y) = \frac{2}{\sigma\sqrt{2\pi}} e^{-\frac{(XY-a)^2}{2\sigma^2}} \cdot x = \frac{2}{\sigma/x\sqrt{2\pi}} e^{-\frac{(Y-a/x)^2}{2\sigma/x^2}}, \qquad (2.15)$$

and the density of the random variable $Y: P(Y) = 0.$, if $Y > a/x = Li(x)/x$.

Let:

$$b = a/x = Li(x)/x, \sigma_1 = \sigma/x = \sqrt{Li(x)/x^2 - Li^2(x)/x^3}. \qquad (2.16)$$

The random variable $Y$ (based on (2.15)) has the same view of the distribution function as the random variable $X_2$, but with the characteristics (2.16). The random variable $Y$ is an analogue of the average density of the number of primes. If you compare the Gaussian density $1/\ln(x)$ with the expectation of a random variable $Y$, which is equal to $b = Li(x)/x$, then based on (2.12) - $1/\ln(x) < Li(x)/x$, for $x < x_0$.

The value of the standard deviation $\sigma_1$ (based on (2.16)) tends rapidly to 0 with increasing $x$. Therefore, the sequence of random variables $Y$ (an analogue of the average density of the number of primes) converges in probability to its mathematical expectation $Li(x)/x$, which is equal to probability of a random integer (not exceed $x$) to be a prime.

Now we consider the random variable $Z = 1/Y$, which is analogous of the mean distance between primes if $x < x_0$, where $x_0$ - the number of Schiusa.

We find the density function of the random variable $Z$, considering that in this case $Z = f(Y)$ is a decreasing function:

$$P_Z(Z) = -P_Y[g(Z)]g'(Z), \qquad (2.17)$$

where $g(Z)$ is the inverse function.

In this case $g(Z) = 1/Z, g'(Z) = -1/Z^2$, therefore based on (2.17) we obtain:

$$P_Z(Z) = \frac{2}{Z^2 \sigma_1 \sqrt{2\pi}} e^{-\frac{(1/Z - a/x)^2}{2\sigma_1^2}} \text{ if } Z \geq c = x/a = x/Li(x), \qquad (2.18)$$

and

$$P_Z(Z) = 0 \text{ If } Z < x/Li(x), \qquad (2.19)$$

where $\sigma_1 = \sigma/x$.

Based on (2.18), (2.19) the expectation of a random variable $Z$ (analogue average distance between prime numbers) is equal to $c = x/Li(x)$.

Based on (2.12) we have:

$$x/Li(x) < x/\pi(x) < \ln(x), \qquad (2.20)$$

where $x/\pi(x)$ is the average distance between primes.

Function $P(Z)$ (if $Z \geq c = x/Li(x)$) reaches a maximum value at the point $c = x/Li(x)$ and based on (2.18), (2.19) is equal to:

$$P(c) = \sqrt{\frac{2/\pi}{1/x^2 Li^3(x) - 1/x^3 Li^2(x)}}. \qquad (2.21)$$

The derivative - $P'(Z) < 0$ (if $Z > c$) and therefore the function $P(Z)$ rapidly decreases to 0.

We use the fact that for the random variable $Y = f(X)$ keep all probability; therefore if the random variable $X$ takes a value $x$ with the probability $p$ that a random variable $Y$ takes the value $y = f(x)$ with the same probability.

Therefore, if the random variable $x_2$ (analogous to the number of primes that do not exceed the value $x < x_0$) is in the range $(Li(x) - \sigma(x), Li(x))$ with a probability of $0,6827$ (for example), then the random variable $Z$ (analogous to the mean distance between the prime numbers not exceeding $x < x_0$) is in the range $Li(x) - \sigma(x), Li(x))$ with the same probability.

Based on the above we consider the following values: the integer value of the standard deviation for the first probabilistic model - $\sigma(x) = [\sqrt{Li(x) - Li^2(x)/x}]$, the expectation of a random variable $Z$ - $x/Li(x)$, the value of the difference between the actual and the calculation values of the mathematical expectation of the average distance between the primes - $x/\pi(x) - x/Li(x)$, the calculated deviation of the random variable $Z$ at the point $Li(x) - \sigma(x)$ from their the expectation - $x/(Li(x) - \sigma(x)) - x/Ln(x)$. These data are presented in Table 2.

Таблица 2

| x | The integer value of the standard deviation for the first probabilistic model | The expectation of the random variable $Z$ | The difference between the actual and the calculated value of the expectation of the random variable $Z$ | The difference between the actual and calculation values of the average distance between primes | Calculated deviation of the random variable $Z$ at the point $Li(x) - \sigma(x)$ from their mathematical expectation |
|---|---|---|---|---|---|
| $10^8$ | 2330 | 17,354 | 17,361 | 0,003 | 0,007 |
| $10^9$ | 7091 | 19,666 | 19,669 | 0,001 | 0,003 |
| $10^{10}$ | 20841 | 21,975 | 21,976 | 0,000 | 0,001 |
| $10^{11}$ | 62856 | 24,283 | 24,284 | 0,000 | 0,001 |

| $10^{12}$ | 612099 | 26, 590 | 26,590 | 0,000 | 0.000 |

We can draw the following conclusions (based on the data of Table 2):

1. The value of difference between the actual and the calculation value of the mathematical expectation of the average distance between primes (for large $x$) is a little.

2. The calculation deviation of the random variable $Z$ at the point $Li(x) - \sigma(x)$ from their mathematical expectation - $x/(Li(x) - \sigma(x)) - x/Ln(x)$ (for large $x$) is a little.

3. The value of the difference between the actual and the calculation values of the mathematical expectation of the average distance between primes does not exceed the deviation of the random variable $Z$ at the point $Li(x) - \sigma(x)$.

## 3. PROBABILISTIC MODELS OF THE DISTRIBUTION OF PRIME TUPLES IN THE NATURAL NUMBERS

The article [5] studied the distribution of the random variable $J(x)$ (analog the number of prime $k$-tuples that do not exceed $x$). Let me remind you that the prime $k$-tuple is a sequence consisting of prime numbers $n, n+2m_1, ... n+2m_{k-1}$, where $m_1 < m_2 < ... < m_{k-1}$. It was also shown that the random variable $J(x)$ has a normal distribution with mean:

$$M_J \approx C(m_1,...m_{k-1}) Li_k(x),$$

where $C(m_1,...,m_{k-1})$ is the constant.

This value $J(x)$ is equal to the number of prime $k$-tuples defined on the basis of Hardy-Littlewoods conjecture. The standard deviation of the random variable $J(x)$ determined by the formula:

$$\sigma_J \approx \sqrt{C(m_1,...,m_{k-1}) Li_k(x) - (C(m_1,...,m_{k-1}))^2 Li_{2k}(x)}.$$

Therefore, the probability density of the distribution of the random variable $J(x)$ is equal to:

$$P_J(J) = \frac{1}{\sigma_J \sqrt{2\pi}} e^{-\frac{(J-M_J)^2}{2\sigma_J^2}}. \tag{3.1}$$

Let us find the probability density of the function of the random variable $G(x) = J(x)/x$. We use the formula:

$$P_G(G) = P_J[G(J)]G'(J), \qquad (3.2)$$

where the random value $G(J)$ is the inverse function.

Based on (22) and (23) we obtain the formula:

$$P_G(G) = \frac{1}{\sigma_J \sqrt{2\pi}} e^{-\frac{(xJ - M_J)^2}{2\sigma_J^2}} \cdot x = \frac{1}{\sigma_J/x \sqrt{2\pi}} e^{-\frac{(J - M_J/x)^2}{2\sigma_J^2/x^2}}. \qquad (3.3)$$

Based on (3.3) the random variable $G(x)$ has a normal distribution as the random variable $J(x)$, but with expectation:

$$M_G = M_J/x \approx \frac{C(m_1, \ldots, m_{k-1}) Li^k(x)}{x} \qquad (3.4)$$

and the standard deviation:

$$\sigma_G \approx \sqrt{C(m_1, \ldots m_{k-1}) Li_k(x)/x^2 - (C(m_1, \ldots m_{k-1}))^2 Li_{2k}(x)/x^3}. \qquad (3.5)$$

The random variable $G(x)$ is the analogue of the average density of the number of prime $k$ - tuples that do not exceed $x$. Based on (3.5) the standard deviation $\sigma_G$ tends rapidly to 0 with increasing $x$. Therefore, the random variable $G(x)$ (analogue of the average density of the number of prime $k$ - tuples that do not exceed $x$) converges in probability to its mathematical expectation $\frac{C(m_1, \ldots, m_{k-1}) Li^k(x)}{x}$, which is equal to the probability of selected at random natural numbers (not exceed $x$): $n, n + 2m_1, \ldots n + 2m_{k-1}$ to be primes

Let me explain the last sentence. It was shown [8] that the density of an integer, positive, strictly increasing sequence (on a finite interval of the natural numbers) is the ultimate probability measure. The sequence of prime $k$ - tuples satisfies the conditions and therefore its density is the ultimate probability measure which is equal to selected at random natural numbers (not exceed $x$): $n, n + 2m_1, \ldots n + 2m_{k-1}$, where $m_1 < m_2 < \ldots < m_{k-1}$, to be primes.

Now we consider the random variable $H = 1/G$, which is analogous to the mean distance between prime $k$ - tuples (not exceed $x$). We find the density of the function of the random variable $H$, considering that $H = f(G)$ is a decreasing function:

$$P_H(H) = -P_G[g(H)]g'(H), \tag{3.6}$$

where $g(H)$ is the inverse function.

We obtain: $g(H) = 1/H, g'(H) = -1/H^2$ in this case. Therefore, based on (3.6) we obtain the formula:

$$P_H(H) = \frac{1}{H^2 \sigma_G \sqrt{2\pi}} e^{-\frac{(1/H - M_G)^2}{2\sigma_G^2}}. \tag{3.7}$$

The expectation of the random variable $H$ (analogue of the mean distance between the prime $k$-tuples not exceed the value $x$) is equal to:

$$M_H = 1/M_G = x/C(m_1, \ldots m_{k-1})Li_k(x). \tag{3.8}$$

The function $P_H(H)$ (3.7) is not symmetrical and does not have a normal distribution. The maximum of function $P_H(H)$ is achieved when $H = M_H$ (3.8) and is equal to:

$$P_H(M_H) = \frac{1}{(M_H)^2 \sigma_G \sqrt{2\pi}} = \frac{(C(m_1, \ldots, m_{k-1})Li_k(x))^2}{x^2 \sigma_G \sqrt{2\pi}}, \tag{3.9}$$

where $\sigma_G$ is defined by the formula (3.5).

If $H$ tends to infinity, then the function $P_H(H)$ tends to 0.

The random variable $H$ (the analog of the average distance between prime $k$-tuples) is a function of the random variable - $J$ (the analog of number of $k$-tuples that do not exceed the value $x$).

We use the fact that for the random variable $Y = f(X)$ all probability values are saved, that if the random variable $X$ takes the value $x$ with the probability $p$ that a random variable $Y$ also takes the value $y = f(x)$ with probability $p$.

For example, let consider the average distance between twin primes, do not exceed $x = 10^6 - 10^8$ (see Table 3). We define the following settings. The value of the standard deviation- $\sigma_J = \sqrt{1,32\ldots Li_2(x) - (1,32\ldots)^2 Li_4(x)}$, the expectation of a random variable $H$ (the analog of the average distance between twin primes) - $M(H) = x/1,32\ldots Li_2(x)$, the value of a random variable $x/(1,32\ldots Li_2(x) + \sigma_J)$, the value of a random variable $x/(1,32\ldots Li_2(x) - \sigma_J)$, the actual number of ordinary twins, does not exceed $x - x/\pi_2(x)$, the actual average distance between twin primes, does not exceed $x - x/\pi_2(x)$, the deviation of the actual average of the distance between twin primes from calculation - $x/\pi_2(x) - M(H)$.

Таблица 3

| x | The value of the standard deviation | The expectation of the random variable $H$ | The value of the random variable $H$ at the point of the expectation plus the standard deviation | The value of the random variable $H$ at the point of the expectation minus the standard deviation | The actual number of twin primes, does not exceed $x$ | $x/\pi_2(x)$ | $x/\pi_2(x) - M(H)$ |
|---|---|---|---|---|---|---|---|
| $10^5$ | 35 | 80,064 | 77,882 | 82,372 | 1224 | 81,699 | 1,635 |
| $10^6$ | 90 | 121,242 | 119,933 | 122,579 | 8169 | 122,414 | 1,1172 |
| $10^7$ | 242 | 170,201 | 169,503 | 170,905 | 58980 | 169,549 | 0,0652 |

We can draw the following conclusions (based on the data of Table 3):

1. The value of the deviation of the actual and calculated average distance between twin primes $x/\pi_2(x) - M(H)$ at large $x$ is a little.

2. The probability that the average distance between twin primes is in the range $(\dfrac{x}{1,32...Li_2(x)+\sigma_J}, \dfrac{x}{1,32...Li_2(x)-\sigma_J})$ is equal to $F(S)$, where $F(S)$ is the value of the standard normal distribution function at the point $S$.

3. The actual average distance between twin primes is in the range $(\dfrac{x}{1,32...Li_2(x)+\sigma_J}, \dfrac{x}{1,32...Li_2(x)-\sigma_J})$ for values $x$ in the range $10^5 - 10^7$.

Now we consider (as an example) the average distance between prime tuples $(p, p+4, p+6)$ (not exceed $x$). Let us define the following settings. The value of the standard deviation - $\sigma_J = \sqrt{C(4,6)Li_3(x) - (C(4,6))^2 Li_6(x)}$, the expectation of a random variable $H$ is the analog of the average distance between prime tuples $(p, p+4, p+6)$: $M(H) = x / C(4,6)Li_3(x)$, the value of a random variable $x/(C(4,6)Li_3(x) + 3\sigma_J)$, the value of a random variable $x/(C(4,6)Li_3(x) - 3\sigma_J)$, the actual number of prime tuples $(p, p+4, p+6)$ (not exceed $x$): $\pi_3(x)$, the actual average distance between prime tuples (do not exceed $x$): $x/\pi_3(x)$, the deviation of the actual average the distance between prime tuples from calculated: $x/\pi_3(x) - M(H)$ (see Table 4).

Таблица 4

| x | The value of the standard deviation | The expectation of the random variable $H$ (an analog of the average distance between prime tuples) | The value of the random variable $H$ at the point of the expectation plus three standard deviations | The value of the random variable $H$ at the point of the expectation minus three standard deviations | The actual number of prime tuples | $x/\pi_3(x)$ | $x/\pi_3(x) - M(H)$ |
|---|---|---|---|---|---|---|---|
| $10^6$ | 16 | 691,563 | 669,344 | 715,308 | 1444 | 692,521 | 0,958 |
| $10^7$ | 38 | 1164,563 | 1164,009 | 1150,086 | 8677 | 1152,472 | -11,536 |
| $10^8$ | 93 | 1802,094 | 1793,079 | 1811,086 | 55556 | 1799,986 | -2,108 |

We can draw the following conclusions (based on the data of Table 4):

1. The value of the deviation of the actual average distance between prime tuples $(p, p+4, p+6)$ from calculated is equal to $|x/\pi_3(x) - M(H)|$ and is relatively small for large $x$.

2. The probability that the average distance between prime tuples $(p, p+4, p+6)$ (for any $x$) to be in range: $(\dfrac{x}{C(4,6)Li_3(x) + S\sigma_J}, \dfrac{x}{C(4,6)Li_3(x) - S\sigma_J})$ is equal to $F(S)$, where $F(S)$ is the value of the function of standard normal distribution at the point $S$.

3. The actual average distance between prime tuples $(p, p+4, p+6)$ is within the range $(\dfrac{x}{C(4,6)Li_3(x) + 3\sigma_J}, \dfrac{x}{C(4,6)Li_3(x) - 3\sigma_J})$ for values $x$ in the range $10^6 - 10^8$.

## 4. CONCLUSION AND SUGGESTIONS FOR FURTHER WORK

The probabilistic estimations can be used in many cases as the only method of analysis of the distribution of primes. Therefore it is opportune to obtain probabilistic estimates of the accuracy of other indicators of the distribution of primes and prime k-tuples.

## 5. ACKNOWLEDGEMENTS

Thanks to everyone who has contributed to the discussion of this paper.